# LARGE DEVIATIONS OF A MODIFIED JACKSON NETWORK: STABILITY AND ROUGH ASYMPTOTICS

By Robert D. Foley[1] and David R. McDonald[2]

*Georgia Institute of Technology and University of Ottawa*

Consider a modified, stable, two node Jackson network where server 2 helps server 1 when server 2 is idle. The probability of a large deviation of the number of customers at node one can be calculated using the flat boundary theory of Schwartz and Weiss [*Large Deviations Performance Analysis* (1994), Chapman and Hall, New York]. Surprisingly, however, these calculations show that the proportion of time spent on the boundary, where server 2 is idle, may be zero. This is in sharp contrast to the unmodified Jackson network which spends a nonzero proportion of time on this boundary.

**1. Introduction.** In this paper we derive the rough (logarithmic) asymptotics for the steady state probability $\pi$ of a particular two node queueing network as the queue at server 1 gets large. The analyzed queueing network is a variation of a two node Jackson queueing network in which server 2 when idle can assist server 1. Allowing one of the servers to help can completely change the behavior of the network. This network clearly exhibits a large deviation phenomenon, which we call a bridge. For certain parameters, as the queue length at node 1 grows, the queue length at node 2 stays small, but generally positive so that server 2 is prevented from helping server 1. Instead of jittering along the $x$-axis, the process skims above the $x$-axis and only rarely touches the axis.

This bridge phenomenon seems to have been somewhat overlooked. In particular, the theory in [8, 11] for analyzing exact asymptotics does not apply. In a companion paper [9], we extend the theory and develop an approach to obtaining the exact asymptotics of networks exhibiting the bridge phenomenon.

Received May 2002; revised February 2004.

[1]Supported in part by the NSF Grant DMI-99-08161.

[2]Supported in part by NSERC Grant A4551.

The bridge phenomenon in the modified Jackson network is not an isolated case—the bridge phenomenon is ubiquitous. Since becoming aware of it, we are encountering it frequently in a variety of contexts. In Section 5 we use this theory to revisit the *bathroom problem* discussed by Shwartz and Weiss [15].

In Section 2 we describe the Jackson network and the modified network. We then discuss the possible large deviation paths for overloading node 1. Section 3 determines the stability conditions of the modified Jackson network. Section 4 analyzes the rough asymptotics of the modified Jackson network; the analysis of the exact asymptotics appears in [9]. Section 5 briefly describes another model where the bridge phenomenon occurs.

## 2. Notation and main results.

Consider a Jackson (1957) network with two nodes. The arrival rate of exogenous customers at nodes 1 and 2 form Poisson processes with rates $\bar{\lambda}_1$ and $\bar{\lambda}_2$, respectively. The service times are independent, exponentially distributed random variables with mean $1/\mu_1$ and $1/\mu_2$, respectively. Each customer's route through the network forms a Markov chain. A customer completing service at node 1 is routed to node 2 with probability $r_{1,2}$ or leaves the system with probability $r_{1,0} := 1 - r_{1,2}$. Routing from node 2 is defined analogously. So without loss of generality, we are assuming $r_{1,1} = r_{2,2} = 0$. The routing process, service processes and arrival processes are independent.

To ensure that the network is open, we assume that $r_{1,2}r_{2,1} < 1$. Since the network is open, the traffic equations

$$(2.1) \qquad \lambda_i = \bar{\lambda}_i + \lambda_{3-i} r_{3-i,i} \qquad \text{for } i = 1, 2,$$

have a unique solution $(\lambda_1, \lambda_2) = ((\bar{\lambda}_1 + \bar{\lambda}_2 r_{2,1})/(1 - r_{1,2}r_{2,1}), (\bar{\lambda}_2 + \bar{\lambda}_1 r_{1,2})/(1 - r_{1,2}r_{2,1}))$. To eliminate degenerate situations, we assume that $\lambda_1 > 0$ and $\lambda_2 > 0$.

The joint queue length process of this Jackson network forms a Markov process with state space $S = \{0, 1, \ldots\}^2$. Define $\rho_i = \lambda_i/\mu_i$, for $i = 1, 2$. From Jackson (1957), it follows that the stationary distribution for the joint queue length process being in the state $(x, y) \in S$ is $(1 - \rho_1)\rho_1^x(1 - \rho_2)\rho_2^y$, provided that the stability conditions $\rho_1 < 1$ and $\rho_2 < 1$ hold.

The network that we analyze is a small change from the above network. Suppose that server 2 has been cross-trained and helps server 1 whenever queue 2 is empty. Let $\mu_1^* \geq \mu_1$ be the combined service effort of the two servers at node 1 when server 2 is empty. The transition rates for the joint queue length process of the modified network are shown in Figure 1. By comparing jump rates, the total number of customers in this modified network is stochastically smaller than the total number in the associated Jackson network. Hence, the modified network will be stable if the associated Jackson network is. However, cross-training server 2 may allow the modified network to be stable even if $\rho_1 > 1$. In particular, if $\rho_2 < 1$ and $\mu_1^* > (\lambda_1 - \mu_1 \rho_2)/(1 - \rho_2)$, then the modified network is stable; see Section 3 for the argument.

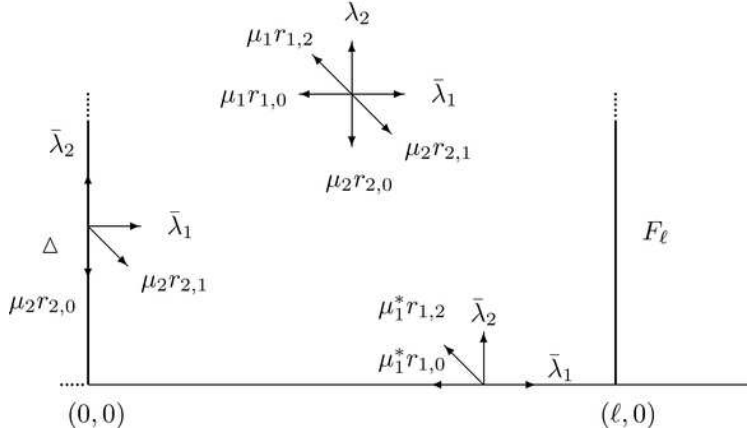

Fig. 1.   *Jump rates for the modified network.*

We are interested in the rare event of a large deviation in the number of customers at node 1; that is, more than $\ell$ customers at node 1 where $\ell$ is large. The steady state probability of this rare event is proportional to the number of visits to $F_\ell \equiv \{(x, y) : x \geq \ell, y \geq 0\}$ between returns to the origin. For the Jackson network, we can determine the most likely path from the origin to $F_\ell$ by looking at the reversed process starting from steady state in $F_\ell$ and look at sample paths that leave $F_\ell$ on the first step and never return. The reversed process has external arrivals entering node $i$ with rate $\lambda_i r_{i,0}$. The service rates at the nodes are unchanged, but the routing probabilities for the reversed process are $r_{j,i}^* = \lambda_i r_{i,j}/\lambda_j$ for $i$ and $j$ in $\{1, 2\}$. Start the reversed process in state $(\ell, y)$, where $\ell$ is large but $y$ is small. As long as there are customers at node 1, customers leave node 1 at rate $\mu_1$. Thus, customers enter node 2 at rate $\lambda_2 r_{2,0} + \mu_1 r_{1,2}^*$. If this rate is less than the maximum rate at which customers can leave node 2, that is,

$$(2.2) \quad \lambda_2 r_{2,0} + \mu_1 r_{1,2}^* < \mu_2 \quad \text{or, equivalently,} \quad \rho_2^{-1} > r_{2,0} + r_{2,1} \rho_1^{-1},$$

then the number of customers at node 2 remains small, and the reversed process starting from $F_\ell$ bounces along the $x$-axis to $(0, 0)$. If the inequality is reversed,

$$\rho_2^{-1} < r_{2,0} + r_{2,1} \rho_1^{-1},$$

then the process starting from $(\ell, 0)$ leaves the $x$-axis and heads roughly northwest (with an easily determined slope) as the customers in node 2 grow until node 1 empties, that is, hits the $y$-axis.

From there, the process bounces along the $y$-axis south to the origin because customers in the time reversed network enter node 1 at rate $\lambda_1 r_{1,0} + \mu_2 r_{2,1}^*$. This rate is less than the maximum rate at which customers can leave node 1, that is,

$$(2.3) \quad \lambda_1 r_{1,0} + \mu_2 r_{2,1}^* < \mu_1 \quad \text{or, equivalently,} \quad \rho_1^{-1} > r_{1,0} + r_{1,2}\rho_2^{-1},$$

given (2.2) fails (otherwise just add the two inequalities together and derive a contradiction).

If the inequality in (2.2) is changed to an equality, then the number in node 2 in the reversed process behaves like a simple, symmetric random walk, which would hit the $y$-axis at a height proportional to $\sqrt{\ell}$. Thus, in the Jackson network, there are three possibilities for the most likely approach from the origin to $F_\ell$, though the approach corresponding to equality in (2.2) occurs only for a set of parameters with Lebesgue measure zero.

Now consider the modified network with $\mu_1^* = \mu_1$. Of course, this is identical to the Jackson network. Suppose the most likely path for a large deviation at node 1 of this Jackson network bounces along the $x$-axis. Let $\mu_1^*$ increase. As $\mu_1^*$ increases, the approach going out the $x$-axis becomes more difficult and may eventually become more difficult than some other approach. In addition to the obvious possibility of going up the $y$-axis, it turns out that there is a third possibility hinted at by the case of equality in (2.2): the process travels along the $x$-axis, but instead of jittering along the $x$-axis, the process skims above and only rarely touches the $x$-axis. This third approach we refer to as a *bridge path*, which is slightly optimistic since we hope to prove properties in a later paper that would justify the word "bridge."

Section 4 contains our preliminary investigation of this modified Jackson network by looking at the behavior under the fluid scaling. For the fluid scaling, speed up the transition rates by $\ell$ and measure customers in units of $1/\ell$, which results in a functional s.l.l.n. In particular, we apply the flat boundary theory of Schwartz and Weiss [15] to obtain rough asymptotics, as well as the fluid scaled large deviation path. In the fluid scaling, both the bridge path and the path that jitters along the $x$-axis collapse to a constant speed line along the $x$-axis, which suggests that the flat boundary theory might not be able to distinguish between the two. However, the calculations also give the proportion of time spent on the boundary. In some cases, the proportion of time spent on the boundary is zero, proving the existence of this third possible approach to $F_\ell$. In fact, we define a *bridge path* to be such a large deviation path which follows a line, for example, an axis though the proportion of time the process spends on the line is zero. Although the term is defined with respect to the fluid scaling, the basis for the term is the conjectured behavior of the unscaled process. Even though we suspect that it has a bridge shape, there are other possibilities. For example, the most likely path when equality holds in (2.2) jitters up the $y$-axis proportional to $\sqrt{\ell}$ before drifting to $F_\ell$ also spends zero time on the boundary and collapses to the $x$-axis under the fluid scaling. We intend to sort out these questions in a future paper.

A *jitter path* follows a line, for example, one of the axes, while spending a nonzero proportion of time on the line. Though jitter path is defined with respect to the fluid scaling, the term reflects the behavior of the unscaled process, which jitters along the line as it travels to $F_\ell$.

We will use the phrase "with large deviation rate $\theta$" to mean that

$$\lim_{\ell \to \infty} \frac{1}{\ell} \log P(W \in F_\ell | W(0) = (0,0)) = -\theta,$$

where $F(\ell)$ is the set of cadlag paths in $S$ starting at the origin and associated with a large deviation of $W$ to $F_\ell$ before returning to the origin; that is, to describe the rough asymptotics. Note $F(\ell)$ is a set of paths hitting the set $F_\ell$.

Basically, we will show that the rough asymptotics can be determined from three points; see Figure 2. The coordinates of the easternmost point of the egg-shaped curve is labelled $\theta^b$. If the curve $M^- = 0$ intersects the egg $M^+ = 0$ between $\theta^b$ and the $\theta_1$-axis, then the intersection is labelled $\theta^j$; otherwise, $\theta^j = \theta^b$. If the horizontal line at height $\log(\rho_2^{-1})$ intersects the egg between $\theta^b$ and the $\theta_1$-axis and (2.3) holds, then the intersection is labelled $\theta^c$; otherwise, $\theta^c = \theta^b$.

The first coordinate of $\theta^b$, $\theta^j$ and $\theta^c$ gives the large deviation rate of the best bridge path, jitter path and cascade paths, respectively. The minimum of the three first coordinates is the rate associated with a large deviation at node 1. Theorem 4 summarizes these results.

**3. A bound and stability.** We will need the following bound in [9]. Since the stability argument and the derivation of the bound use the same coupling, we have included both in this section.

LEMMA 1. *For the stable, modified network,*

$$(3.1) \qquad \sum_{j \geq y} \pi(0, j) \leq c \rho_2^y.$$

PROOF. The proof is divided into two parts. First, we use a coupling argument to consider the case $\lambda_1 < \mu_1$. We associate a region $R(x, y)$ with each point $(x, y) \in S$. The argument shows that there is a coupling of the queue length processes of the Jackson and modified networks so that if the modified network is in state $(x, y)$, then the state of the Jackson network is in $R(x, y)$. It immediately follows that the stationary probability of the modified network being in state $(x, y)$ is bounded by the stationary probability that the Jackson network is in $R(x, y)$. Similarly, it follows that the stationary probability that the modified network is in $B \subset S$ is less than or equal to the stationary probability that the Jackson network is in $\bigcup_{(x,y)\in B} R(x, y)$.

Now we describe the coupling. We can consider them to be a pair of discrete time Markov chains, $W[n]$ for the Jackson network and $Y[n]$ for the modified network, subordinated to a common Poisson process with rate $\lambda_1 + \lambda_2 + \mu_1^* + \mu_2$, which we assume without loss of generality to be one. Basically, $W[n]$ and $Y[n]$ "attempt" to move in the same direction. More precisely, generate an i.i.d. sequence of random variables (directions) taking values $\{E, N, W, NW, S, SE, W^*, NW^*\}$ with probabilities $\{\lambda_1, \lambda_2, \mu_1(1 - r_{1,2}), \mu_1 r_{1,2}, \mu_2(1 - r_{2,1}), \mu_2 r_{2,1}, (\mu_1^* - \mu_1)(1 - r_{1,2}), (\mu_1^* - \mu_1) r_{1,2}\}$, respectively.

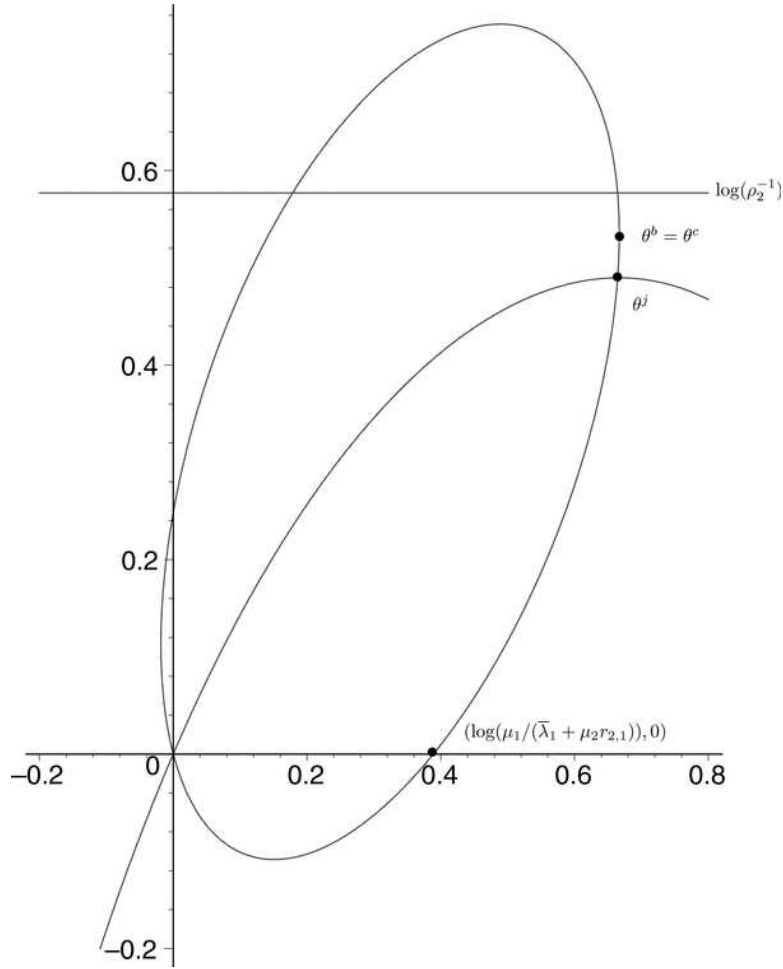

Fig. 2.   $M^+ = 0$ is the egg-shaped curve; the other curve is $M^- = 0$.

For the $n$th step, both $W[n]$ and $Y[n]$ "attempt" to take a single step in the direction given by the $n$th random direction. The directions $W^*$ and $NW^*$ indicate west and northwest, but only for the modified network when it is on the $x$-axis; otherwise, $W^*$ and $NW^*$ indicate that the process stays put. By "attempt," we mean that the process moves to the neighboring state in that direction unless it would result in the process leaving the state space; for example, if either process were on the $y$-axis and the attempted step were to be in the direction $W$. Note that if both processes move, they move in the same direction. However, one process may move while the other stays put.

For the regions, define

(3.2)
$$\begin{aligned} R(x,y) = &\{(i,j) \in S | i \geq x, j \geq (y-1)^+\} \\ &\cup \{(i,j) \in S | (i,j) = (x-1,y)\}. \end{aligned}$$

Note that for $(x,y) \geq (1,1)$, both $(x-1,y)$ and $(x,y-1)$ are in $R(x,y)$, but $(x-1,y+1)$ and $(x-1,y-1)$ are not. We claim that if $W[n] \in R(Y[n])$, then $W[n+1] \in R(Y[n+1])$. This is clear when both processes are in the interior since they both move in the same direction. Consider the case $Y[n] = (1,y)$ with $y > 1$ and $W[n] = (0,y)$. Both processes move in the same direction unless the direction is $W$ or $NW$; in all cases, $W[n+1]$ will still be in $R(Y[n+1])$. Now consider the case when $Y[n] = (1,0)$ and $W[n] = (1,0)$. Consider the movements $NW^*$, then $W$ and then $SE$. This makes $W[n+3] = (0,0)$ and $Y[n+3] = 1$ and $W[n+k]$ will still be in $R(Y[n+k])$ for $k = 1,2,3$. This trajectory explains why $R(x,y)$ is defined as it is. We leave it to the reader to finish checking the claim.

Now start $W[0]$ and $Y[0]$ off in the same state with distribution $\pi$. Notice that

$$\begin{aligned} \sum_{i \geq 0, j \geq y} \pi(i,j) &= P(Y[n] \in R(0,y+1)) \\ &= P(W[n] \in R(Y[n]), Y[n] \in R(0,y+1)) \\ &\leq P(W[n] \in R(0,y)) < \rho_2^{(y-1)^+} \end{aligned}$$

for sufficiently large $n$ because $W[n]$ converge to steady state. This completes the argument under the condition $\lambda_1 < \mu_1$. [$R(x,y)$ is the smallest set that will work for this coupling. This can be seen by starting both chains at the origin and arguing that there is a sequence of random directions such that $(Y[n], W[n]) = (x,y,i,j)$ for every $(x,y) \in S$ and $(i,j) \in R(x,y)$. Most points can be reached by having both processes move sufficiently far east, then the modified network moves west back to the origin, then the Jackson network moves sufficiently far northwest, and finally both processes move north and east sufficiently far. The few remaining points can be reached by getting the Jackson network to the origin while the modified network is at $(0,1)$.]

Now consider the case when $\lambda_1 \geq \mu_1$. The number of customers at node 2 is stochastically smaller than a birth–death process on the nonnegative integers with birth rate $\bar{\lambda}_2 + \mu_1 r_{1,2}$ and death rate $\mu_2$ on the positive integers, and birth rate $\bar{\lambda}_2 + \mu_1^* r_{1,2}$ in state 0. Thus, the probability of more than $y$ customers at node 2 is smaller than $c[(\bar{\lambda}_2 + \mu_1 r_{1,2})/\mu_2]^y$ for a suitably chosen constant $c$. Using (2.1) and that $\lambda_1 \geq \mu_1$, it follows that $[(\bar{\lambda}_2 + \mu_1 r_{1,2})/\mu_2] \leq \rho_2$. □

### 3.1. *Stability of the modified network.*

PROPOSITION 1. *The joint queue length process of the modified network is positive recurrent if $\lambda_2 < \mu_2$ and $\lambda_1 < \rho_2 \mu_1 + (1-\rho_2)\mu_1^*$. If either inequality is reversed, then the process is transient.*

PROOF. The result follows from comparing the modified and Jackson networks in the case $\rho_1 < 1$; for the coupling, see the first part of the proof of Lemma 1. Now consider the case where $\lambda_1 \geq \mu_1$. The following definitions of $T_n$ and $X_n$ are only used in this proof. Let $T_n$ denote the time that the $n$th busy period starts at node 2, and let $X_n$ be the number of customers in queue 1 just prior to the start of the $n$th busy period at node 2. The process $X_0, X_1, \ldots$ is a Markov chain. To prove the result, it suffices to show that $\mathrm{E}[T_{n+1} - T_n | X_n] < \infty$ and that $X_0, X_1, \ldots$ is positive recurrent.

The random variable $T_{n+1} - T_n$ represents the $n$th busy cycle at node 2, which is the sum of the $n$th busy period and the $n$th idle period at node 2. Since $\lambda_2 > 0$, the expected length of the $n$th idle period is finite. The length of the $n$th busy period at node 2 is stochastically increasing in $X_n$; hence, it is stochastically smaller than the busy period in an $M/M/1$ queue with service rate $\mu_2$ and arrival rate $\bar{\lambda}_2 + \mu_1 r_{1,2}$. Fortunately, $\bar{\lambda}_2 + \mu_1 r_{1,2} < \mu_2$ since $\bar{\lambda}_2 + \mu_1 r_{1,2} < \lambda_2$ when $\lambda_1 > \mu_1$. The expected length of this bounding busy period is $1/(\mu_2 - (\bar{\lambda}_2 + \mu_1 r_{1,2}))$.

To show that $X_0, X_1, \ldots$ is positive recurrent, we will show that for all $X_n$ greater than some large constant $M$, $\mathrm{E}[X_{n+1} - X_n | X_n] < -\varepsilon < 0$, which is a Foster–Lyapunov type condition guaranteeing stability. The constant $M$ can be chosen large enough so that the busy period at node 2 is arbitrarily close to the busy period of an $M/M/1$ queue with service rate $\mu_2$ and arrival rate $\bar{\lambda}_2 + \mu_1 r_{1,2}$ and that the departure process from node 1 during this busy period is arbitrarily close to a Poisson process with rate $\mu_1$. The change in the queue length at node 1, $X_{n+1} - X_n$, can be decomposed into the sum of the change during the busy period at 2, and the change during the idle period at 2. Hence, for $X_n > M$, $\mathrm{E}[X_{n+1} - X_n | X_n]$ is arbitrarily close to

$$(3.3) \qquad \frac{(\bar{\lambda}_1 + \mu_2 r_{2,1} - \mu_1)}{(\mu_2 - (\bar{\lambda}_2 + \mu_1 r_{1,2}))} + \frac{(\bar{\lambda}_1 - \mu_1^*)}{(\bar{\lambda}_2 + \mu_1^* r_{1,2})},$$

which is strictly less than zero if $\mu_1^* > (\lambda_1 - \mu_1 \rho_2)/(1 - \rho_2)$ and is equivalent to $\rho_2 \mu_1 + (1 - \rho_2)\mu_1^* > \lambda_1$.

To show transience, first assume that $\lambda_2 < \mu_2$, but $\lambda_1 > \rho_2 \mu_1 + (1 - \rho_2)\mu_1^*$. From the argument above, for $X_n \geq M$ sufficiently large, $\mathrm{E}[X_{n+1} - X_n | X_n] > \varepsilon > 0$; hence, for $M$ sufficiently large, $X_n$ behaves arbitrarily closely to a random walk with strictly positive drift and is transient.

Now consider the case when $\lambda_2 > \mu_2$. Note that $\lambda_i$ represents the long run average arrival and departure rate from node $i$ assuming that the network, either Jackson or modified, is recurrent. Assume $X_n$ is recurrent; otherwise, we would be done. Since $X_n$ is recurrent, the departure rate from node 1 equals the arrival rate. If this rate is less than $\lambda_1$, then the queue length at node 2 must be diverging. However, if the departure rate from node 1 is $\lambda_1$, then node 2 is also diverging since the arrival rate would be $\lambda_2$. In either case, the network is transient. □

## 4. Flat boundary approach.

We are interested in describing how a modified Jackson network overloads. In particular, we will be interested in the rare event when the system starts out empty and node 1 reaches a level $\ell$ before the system empties again.

If we speed up the jump rates by a factor $\ell$, but reduce the jumps by a factor $1/\ell$, then we get the scaled process $W_\ell$. From the theory in [15],

$$\lim_{\ell \to \infty} \frac{1}{\ell} \log P(W \in F(\ell) | W(0) = (0,0))$$

$$= \lim_{\ell \to \infty} \frac{1}{\ell} \log P(W_\ell \in F(1) | W_\ell(0) = (0,0))$$

(4.1)
$$= -\inf_{p \in F}[I(p)],$$

where $I(p) = \int_0^T \Lambda(\frac{dp}{ds}(s), p(s)) \, ds$ is the rate function associated with a path $p$ in the set of absolutely continuous paths $F$ starting from $(0,0)$ which hits the set $\{(1, y) : y \geq 0\}$ at some time $T$ before returning to $(0,0)$. Note that after time $T$, the path follows the natural drift path back to the origin, and $\Lambda(\frac{dp}{ds}(s), p(s)) = 0$ along the drift path. We will use better, cheaper and smaller action, synonymously.

Since $W$ has constant (but different) jump rates on and off the flat boundary, the local rate function $\Lambda(\vec{v}, w) = \Lambda^+(\vec{v})$ if $w$ is in the interior; that is, if $w = (x, y)$ with $y > 0$ and $\Lambda(\vec{v}, w) = \Lambda^-(\vec{v})$ if $w = (x, 0)$. It follows from the calculus of variations that the cheapest path in $F$ is a sequence of line segments of constant speed which changes direction only on the $y$-axis. [If a path changes direction on the $x$-axis, either leaving the interior to travel along the $x$-axis or by leaving the $x$-axis at some point other than the origin to travel through the interior, there is a cheaper path from the origin to $\{(1, y) : y \geq 0\}$; the cheaper path might hit $\{(1, y) : y \geq 0\}$ at a different point like $(1, 0)$.] Consequently, the cheapest path in $F$ must lie in one of the following sets $F_i, F_x, F_c$, where:

$F_i$ is the set of all constant speed paths with positive slope across the *interior* until hitting $(1, y)$ with $y > 0$,

$F_x$ is a constant speed path jittering or forming a bridge along the $x$-axis until hitting $(1, 0)$, and

$F_c$ is a constant speed path jittering up the $y$-axis, then changes direction and heads for the point $(1, 0)$ at constant velocity. Thus, the customers first build up in node 2, and then *cascade* into 1.

In each of these three cases, we will be able to reduce the problem to a differentiable, constrained nonlinear optimization problem. We use the Karush–Kuhn–Tucker conditions, which are given in a variety of texts including [3, 13], to determine the minimal action in each of the three cases. These results are then combined to determine (4.1).

### 4.1. A bridge is better than any path through the interior.

In this section we consider the interior paths $F_i \subset F$. That is, we consider paths that initially have the form $p(s) = (v_1 s, v_2 s)$ until hitting $(1, v_2/v_1)$ at time $T = 1/v_1$, where $v_1 > 0$ and $v_2 > 0$. After time $T$, the path follows the natural drift path until reaching the origin. We will find the $\inf_{p \in F_i}[I(p)]$ and show that this inf is not attained in $F_i$. However, there will be a bridge path which attains the inf. Hence, there will always be a bridge that is better than every interior path.

Define the log moment generating function of the compound Poisson process associated with jumps in the interior

$$M^+(\theta_1, \theta_2) = \bar{\lambda}_1(e^{\theta_1} - 1) + \mu_1 r_{1,0}(e^{-\theta_1} - 1) + \mu_1 r_{1,2}(e^{\theta_2 - \theta_1} - 1)$$
$$+ \bar{\lambda}_2(e^{\theta_2} - 1) + \mu_2 r_{2,0}(e^{-\theta_2} - 1) + \mu_2 r_{2,1}(e^{\theta_1 - \theta_2} - 1).$$

The Hessian of $M^+$ is positive definite; hence, $M^+$ is strictly convex. We also know that $M^+(0, 0) = 0$. Now, we argue that there exists a point $(\hat{\theta}_1, \hat{\theta}_2) > (0, 0)$ with $M^+(\hat{\theta}_1, \hat{\theta}_2) < 0$. To see this, first consider the case $\rho_1 < 1$. In this case, $\nabla M^+(0, 0) \cdot (1, 1) < 0$. Hence, $(1, 1)$ is a decreasing direction. Now assume $\rho_1 \geq 1$. In this case, $(0, 1)$ and $(\varepsilon, 1)$ are decreasing directions for a suitably small $\varepsilon > 0$.

Next, the local rate function in the interior in direction $\vec{v} = (v_1, v_2)$ is

$$\Lambda^+(\vec{v}) = \sup_{\theta_1, \theta_2}(\theta_1 v_1 + \theta_2 v_2 - M^+(\theta_1, \theta_2)),$$

which is clearly convex. For a path $p \in F_i$ with velocity $\vec{v}$, $\Lambda(\frac{dp}{ds}(s), p(s)) = \Lambda^+(\vec{v})$ for $s > 0$. Hence, $I(p) = \Lambda^+(\vec{v})/v_1$ is a good rate function; see Theorem 5.1 of [15]. The remainder of this section is devoted to finding $\inf_{p \in F_i}[I(p)] =$ $\inf_{v_1>0,v_2>0} \Lambda^+(v_1,v_2)/v_1$. The argument consists of 3 steps. First, we argue that $\inf_{v_1>0,v_2>0} \Lambda^+(v_1,v_2)/v_1 = \min_{v_1>0,v_2\geq 0} \Lambda^+(v_1,v_2)/v_1$. Second, we argue that every local minimum is a KKT (Karush–Kuhn–Tucker) point, which will be defined shortly. Third, we argue that there is exactly one KKT point $(v_1,v_2) \in (0,\infty) \times [0,\infty)$; hence, this KKT point must be the global minimum of $\min_{v_1>0,v_2\geq 0} \Lambda^+(v_1,v_2)/v_1$.

For $v_1 \geq 0$ and $v_2 \geq 0$, clearly, $\Lambda^+(v_1,v_2) \geq -M^+(\hat{\theta}_1,\hat{\theta}_2) > 0$. Hence, $\Lambda^+(v_1,v_2)/v_1$ is bounded below by 0 and explodes as $v_1$ decreases to zero. Furthermore, by Proposition 3.1 in [4], $\Lambda^+(v_1,v_2)/v_1$ explodes as the norm of $\vec{v}$ becomes large. Also, for each $v_1 > 0$, we know that $\Lambda^+(v_1,v_2)/v_1$ is continuous in $v_2$ as $v_2$ converges to 0. Hence, $\inf_{v_1>0,v_2>0} \Lambda^+(v_1,v_2)/v_1$ must be a local minimum (with $v_1 > 0$) of

$$(4.2) \qquad\qquad \min \Lambda^+(v_1,v_2)/v_1$$

$$(4.3) \qquad\qquad \text{s.t.} \quad v_2 \geq 0.$$

It is known that every locally optimal solution to a constrained, differentiable, nonlinear optimization problem with linear constraints must be a KKT point; see 14.37 of [13] or Section 3.5 of [3]. Now, we argue that our constrained, nonlinear optimization problem is differentiable.

The compound Poisson distribution has an infinite support over the integers; hence, by Proposition 3.1 in [4], there are unique values $\theta_1^+(\vec{v}), \theta_2^+(\vec{v})$, such that

$$(4.4) \qquad \Lambda^+(\vec{v}) = \theta_1^+(\vec{v})v_1 + \theta_2^+(\vec{v})v_2 - M^+(\theta^+(\vec{v})).$$

To explore the relationship between $\vec{v}$ and $\theta^+(\vec{v})$, fix $\vec{v}$ and find $\theta^+$ that maximizes $\theta_1 v_1 + \theta_2 v_2 - M^+(\theta_1,\theta_2)$. Taking derivatives and setting them equal to zero yields

$$(4.5) \quad v_1 = \frac{\partial M^+(\theta^+)}{\partial \theta_1} = (\bar{\lambda}_1 e^{\theta_1^+} - \mu_1 r_{1,0} e^{-\theta_1^+} - \mu_1 r_{1,2} e^{\theta_2^+-\theta_1^+} + \mu_2 r_{2,1} e^{\theta_1^+-\theta_2^+}),$$

$$(4.6) \quad v_2 = \frac{\partial M^+(\theta^+)}{\partial \theta_2} = (\bar{\lambda}_2 e^{\theta_2^+} - \mu_2 r_{2,0} e^{-\theta_2^+} + \mu_1 r_{1,2} e^{\theta_2^+-\theta_1^+} - \mu_2 r_{2,1} e^{\theta_1^+-\theta_2^+}).$$

Thus, $\theta_1^+(v_1,v_2)$ and $\theta_2^+(v_1,v_2)$ determine $v_1$ and $v_2$. Furthermore, this mapping from $\theta^+$ to $\vec{v}$ is a smooth bijection. Since the Jacobian, which is positive definite since it is also the Hessian of the strictly convex $M^+$, has a nonzero determinant, it follows from the inversion theorem (see [2]), that $\theta_1^+(v_1,v_2)$ and $\theta_2^+(v_1,v_2)$ are smooth functions of $(v_1,v_2)$. Hence, our nonlinear programming problem is differentiable.

For our problem, $(v_1,v_2)$ is a KKT point if there exists a Lagrange multiplier $u$ with

$$(4.7) \qquad uv_2 = 0 \qquad\qquad \text{(complementary slackness)},$$

$$(4.8) \qquad u \geq 0 \qquad\qquad \text{(sign restriction)},$$

$$(4.9) \qquad \begin{pmatrix} 0 \\ u \end{pmatrix} = \nabla(\Lambda^+(v_1, v_2)/v_1) \qquad \text{(gradient equation)},$$

$$(4.10) \qquad v_2 \geq 0 \qquad\qquad \text{(constraint)}.$$

Among other things, the proof of our next result shows that any KKT point $(v_1, v_2) \in (0, \infty) \times [0, \infty)$ must have $u > 0$; hence, complementary slackness implies that $v_2 = 0$. Consequently, there is no "best" path in $F_i$; that is, the $\inf_{p \in F_i} I(p)$ is not attained in $F_i$.

THEOREM 1.  *For the paths in the interior, we have*

$$\inf_{p \in F_i}[I(p)] = \inf_{v_1 > 0, v_2 > 0} \Lambda^+(v_1, v_2)/v_1 = \theta_1^b,$$

*where* $\theta^b = \theta^+(\vec{v}^b)$ *and* $\vec{v}^b \equiv (v_1^b, v_2^b)$ *is the unique KKT point in* $(0, \infty) \times [0, \infty)$ *point for* (4.2) *and* (4.3)*. The point* $\vec{v}^b$ *is the unique solution to*

$$(4.11) \qquad\qquad v_1^b > 0,$$

$$(4.12) \qquad\qquad v_2^b = 0,$$

$$(4.13) \qquad\qquad M^+(\theta^+(\vec{v}^b)) = 0.$$

*Equivalently, but more usefully,* $\theta^b$ *is the unique solution to*

$$(4.14) \qquad\qquad \theta_1^b > 0,$$

$$(4.15) \qquad\qquad \frac{\partial[M^+(\theta_1^b, \theta_2^b)]}{\partial \theta_2} = 0,$$

$$(4.16) \qquad\qquad M^+(\theta^b) = 0.$$

PROOF.  First,

$$(4.17) \qquad\qquad \frac{\partial[\Lambda^+(v_1, v_2)]}{\partial v_i} = \theta_i^+(v_1, v_2),$$

which can be seen by starting with (4.4) and using the left most equation in (4.5) and (4.6). It follows that

$$(4.18) \qquad \frac{\partial[\Lambda^+(v_1, v_2)/v_1]}{\partial v_1} = [v_1 \theta_1^+(v_1, v_2) - \Lambda^+(v_1, v_2)]/v_1^2$$

and

$$(4.19) \qquad \frac{\partial[\Lambda^+(v_1, v_2)/v_1]}{\partial v_2} = \theta_2^+(v_1, v_2).$$

In the remainder of this proof assume that $(v_1, v_2)$ is a KKT point, and we know that there exists at least one such point. Thus, using (4.19), the

Lagrange multiplier is $u = \theta_2^+(v_1, v_2)$, and using complementary slackness, (4.18) becomes

$$\frac{\partial[\Lambda^+(v_1, v_2)/v_1]}{\partial v_1} = M^+(\theta^+(\vec{v}))/v_1^2.$$

Now, we argue that $u > 0$. Assume the contrary, that is, that $u = 0$. Thus, the gradient of $\Lambda^+(v_1, v_2)/v_1$ at $(v_1, v_2)$ is zero, which means that $M^+(\theta^+(\vec{v})) = 0$ and $\theta_2^+(v_1, v_2) = 0$. Solving $M^+(\theta_1^+, 0) = 0$ yields $\exp(\theta_1^+) = \mu_1/(\bar{\lambda}_1 + \mu_2 r_{2,1})$.

Substituting into (4.5) and (4.6) yields

$$v_1(\theta_1^+, 0) = \mu_1 - (\bar{\lambda}_1 + \mu_2 r_{2,1})$$

and

$$
\begin{aligned}
v_2(\theta_1^+, 0) &= \bar{\lambda}_2 - \mu_2 r_{2,0} + \mu_1 r_{1,2} \frac{(\bar{\lambda}_1 + \mu_2 r_{2,1})}{\mu_1} - \mu_2 r_{2,1} \frac{\mu_1}{(\bar{\lambda}_1 + \mu_2 r_{2,1})} \\
&= \bar{\lambda}_2 - \mu_2 r_{2,0} + r_{1,2} \bar{\lambda}_1 + \mu_2 r_{2,1} r_{1,2} - \mu_2 r_{2,1} \frac{\mu_1}{(\bar{\lambda}_1 + \mu_2 r_{2,1})} \\
&= (\bar{\lambda}_2 + \bar{\lambda}_1 r_{1,2}) - \mu_2(1 - r_{2,1}) + \mu_2 r_{2,1} r_{1,2} - \mu_2 r_{2,1} \frac{\mu_1}{(\bar{\lambda}_1 + \mu_2 r_{2,1})} \\
&= \lambda_2(1 - r_{1,2} r_{2,1}) - \mu_2(1 - r_{1,2} r_{2,1}) + \mu_2 r_{2,1} \left(1 - \frac{\mu_1}{(\bar{\lambda}_1 + \mu_2 r_{2,1})}\right) \\
&= (\lambda_2 - \mu_2)(1 - r_{1,2} r_{2,1}) + \mu_2 r_{2,1} \left(1 - \frac{\mu_1}{(\bar{\lambda}_1 + \mu_2 r_{2,1})}\right) \\
&= (\lambda_2 - \mu_2)(1 - r_{1,2} r_{2,1}) - \mu_2 r_{2,1} \frac{v_1(\theta_1^+, 0)}{(\bar{\lambda}_1 + \mu_2 r_{2,1})} \\
&< 0 \qquad \text{since } \lambda_2 < \mu_2, v_1 > 0,
\end{aligned}
$$

which is a contradiction. Hence, $u > 0$, and from complementary slackness, we have $v_2 = 0$.

Equations (4.12) and (4.13) follow from the gradient equation, and are rephrased in terms of $\theta^b$ in (4.15) and (4.16). Equation (4.14) follows from (4.11) since (4.18) must equal 0. Thus, the two sets of equations are equivalent. To show uniqueness, consider the second set of equations and recall that $M^+$ is a strictly convex function with $M^+(0,0) = 0$. The set of all $(\theta_1, \theta_2)$ such that $M^+(\theta_1, \theta_2) = 0$ is the boundary of a strictly convex set, the egg shaped region in Figure 2, containing $(0,0)$ and $(\hat{\theta}_1, \hat{\theta}_2) > (0,0)$. There are exactly two points on the boundary of this convex set that are tangent to vertical lines, that is, satisfy (4.15), but only one of the two, the eastern most point on the boundary, satisfies (4.14). $\quad\square$

4.2. *The jitter path on the x-axis* In this section we consider paths that bounce along the $x$-axis; that is, we consider paths in $F$ that initally have the form $p(s) = (v_1 s, 0)$ until hitting $(1, 0)$ at time $T = 1/v_1$, where $v_1 > 0$ and after time $T$, the path follows the natural drift path until reaching the origin. Let $F_x \subset F$ denote the set of all such paths, which will be the jitter paths and the bridge path. To analyze these paths, we view the $x$-axis as a flat boundary as in Definition 8.7 in [15] and $W$ as a flat boundary process.

Define the log moment generating functions of the compound Poisson process associated with jumps on the $x$-axis,

$$M^-(\theta_1, \theta_2) = \bar{\lambda}_1(e^{\theta_1} - 1) + \bar{\lambda}_2(e^{\theta_2} - 1) + \mu_1^* r_{1,0}(e^{-\theta_1} - 1) + \mu_1^* r_{1,2}(e^{\theta_2 - \theta_1} - 1).$$

The associated local rate function is

$$\Lambda^-(v_1, v_2) = \sup_{\theta_1, \theta_2}(\theta_1 v_1 + \theta_2 v_2 - M^-(\theta_1, \theta_2)).$$

Using the same arguments as in the previous section, there exists a unique pair $\theta^-(\vec{v})$ such that

$$\Lambda^-(v_1, v_2) = \theta_1^-(v_1, v_2)v_1 + \theta_2^-(v_1, v_2)v_2 - M^-(\theta^-(\vec{v})).$$

The local rate function for the path with velocity $v = (v_1, 0)$ is given by

$$\Lambda^*(v_1, 0) = \inf_{0 \leq \beta \leq 1, \beta \vec{v}^+ + (1-\beta)\vec{v}^- = (v_1, 0)}(\beta \Lambda^+(\vec{v}^+) + (1-\beta)\Lambda^-(\vec{v}^-)),$$

which is a good rate function; see (v) of Lemma 8.20 of [15]. Intuitively, the path is a mixture with $\beta$ representing the proportion of time above the $x$-axis, while $1 - \beta$ is the proportion of time on the $x$-axis. The bridge path has $\beta = 1$; jitter paths have $\beta < 1$. Finally, we must calculate $\inf_{p \in F_x}[I(p)]$; that is,

$$(4.20) \qquad \inf_{v_1 > 0, 0 \leq \beta \leq 1, \beta \vec{v}^+ + (1-\beta)\vec{v}^- = (v_1, 0)} f(\beta, v_1, v_1^+, v_2^+, v_1^-, v_2^-),$$

where $f(\beta, v_1, v_1^+, v_2^+, v_1^-, v_2^-) = (\beta \Lambda^+(\vec{v}^+) + (1-\beta)\Lambda^-(\vec{v}^-))/v_1$. Our argument will be similar to the last section. First, we argue that (4.20) equals

$$(4.21) \qquad \min_{v_1 > 0, 0 < \beta \leq 1, \beta \vec{v}^+ + (1-\beta)\vec{v}^- = (v_1, 0)} f(\beta, v_1, v_1^+, v_2^+, v_1^-, v_2^-).$$

Just as for $\Lambda^+$, $\Lambda^-(\vec{v}^-)$ goes to infinity as $|\vec{v}^-|$ diverges. Furthermore, $\Lambda^-(\vec{v}^-)$ goes to infinity as $v_2^- \downarrow 0$ and is infinite if $v_2^- \leq 0$. Hence, $f$ goes to infinity as $\beta \downarrow 0$. If $\beta > 0$, $f$ goes to infinity as $v_1 \downarrow 0$. Hence, (4.20) must be a local minimum of

$$(4.22) \qquad \min f(\beta, v_1, v_1^+, v_2^+, v_1^-, v_2^-)$$

$$(4.23) \qquad \text{s.t.} \quad g_1(\beta, v_1, v_1^+, v_2^+, v_1^-, v_2^-) \equiv \beta \leq 1,$$

$$(4.24) \qquad \qquad g_2(\beta, v_1, v_1^+, v_2^+, v_1^-, v_2^-) \equiv -v_1 + \beta v_1^+ + (1-\beta)v_1^- = 0,$$

$$(4.25) \qquad \qquad g_3(\beta, v_1, v_1^+, v_2^+, v_1^-, v_2^-) \equiv \beta v_2^+ + (1-\beta)v_2^- = 0.$$

To show that this constrained, nonlinear optimization problem is differentiable is almost identical to the argument in the previous section including the argument that $\theta_i^-(v_1, v_2)$ is smooth. The constraints are no longer linear since there are terms like $\beta v_1^+$. However, the gradients of the three constraint equations are linearly independent; hence, from 14.37 of [13] or Section 3.5 of [3], all points are regular and every local minimum must be a KKT point. The remainder of the argument is to determine the KKT points $(\beta, v_1, v_1^+, v_2^+, v_1^-, v_2^-) \in (0, 1] \times (0, \infty) \times (-\infty, \infty) \times (-\infty, 0] \times (-\infty, \infty) \times [0, \infty)$. Note that we did not include $v_2^+ > 0$ and $v_2^- < 0$; these will be suboptimal since $\Lambda^-(v_1^-, v_2^-) = \infty$ when $v_2^- \leq 0$.

For our problem, $(\beta, v_1, v_1^+, v_2^+, v_1^-, v_2^-)$ is a KKT point if there exists Lagrange multipliers $u_1$, $u_2$ and $u_3$ satisfying

$$u_1(1 - \beta) = 0 \qquad \text{(complementary slackness)},$$

$$u_1 \leq 0 \qquad \text{(sign restriction)},$$

$$\sum_{i=1}^{3} u_i \nabla g_i = \nabla f \qquad \text{(gradient equation)},$$

$$\beta \leq 1 \qquad \text{(constraint 1)},$$

$$-v_1 + \beta v_1^+ + (1 - \beta)v_1^- = 0 \qquad \text{(constraint 2)},$$

$$\beta v_2^+ + (1 - \beta)v_2^- = 0 \qquad \text{(constraint 3)}.$$

Now assume that $(\beta, v_1, v_1^+, v_2^+, v_1^-, v_2^-)$ is a KKT point. The bottom four components of the gradient equation imply that $v_1 u_2 = \theta_1^+(v^+) = \theta_1^-(v^-)$ and $v_1 u_3 = \theta_2^+(v^+) = \theta_2^-(v^-)$. Using these in the second component of the gradient equation implies that $\beta M^+(\theta^+(\vec{v}^+)) + (1 - \beta)M^-(\theta^+(\vec{v}^+)) = 0$. Since $v_1 > 0$, the second component of the gradient equation also implies that $u_2 > 0$, which means that $v_1 u_2 = \theta_1^+(v^+) > 0$. The first component of the gradient equation reduces to $u_1 = -[M^+(\theta^+(\vec{v}^+)) - M^-(\theta^+(\vec{v}^+))]/v_1$.

First, consider the case when $\beta < 1$. By complementary slackness, $u_1 = 0$ implying that

$$(4.26) \qquad M^+(\theta^+(\vec{v}^+)) = M^-(\theta^+(\vec{v}^+)) = 0.$$

These level sets enclose convex regions which intersect at $(0, 0)$ and possibly other points $(\theta_1^+, \theta_2^+)$. However, since $v_2^+ = \frac{\partial M^+(\theta^+)}{\partial \theta_2}$ must be (strictly) negative, we are only interested in solutions to (4.26) in the lower portion of the egg in Figure 2; that is, going clockwise from (but not including) $\theta^b$ along $M^+ = 0$ to the other solution of (4.15), which is the western most point of the egg shaped region defined by $M^+ = 0$. Since $\theta_1^+(v^+) > 0$, we can further restrict the region to the arc going clockwise from the $\theta^b$ to the origin.

By calculation the points [besides $(0,0)$] where $M^+(\theta) = 0$ and $M^-(\theta) = 0$ cut the $\theta_1$-axis are $(\log(\mu_1/(\bar{\lambda}_1 + \mu_2 r_{2,1})),0)$ and $(\log(\mu_1^*/\bar{\lambda}_1),0)$, respectively. Since the first coordinate of the latter is positive and greater than the first coordinate of the former, it follows that $\theta_2^+ > 0$ if $\theta_1^+ > 0$; thus, we can restrict attention to the arc going clockwise from $\theta^b$ to the $\theta_1$-axis at $(\max(\log(\mu_1/(\bar{\lambda}_1 + \mu_2 r_{2,1})),0),0)$. There can be at most one such point since $\frac{\partial M^-(\cdot)}{\partial \theta_2} > 0$. If there is such a point, label this point as $\theta^j = (\theta_1^j, \theta_2^j)$; otherwise, define $\theta^j = \theta^b$. Note that there can be no other KKT points with $\beta < 1$. To see that $\theta^j$ determines a KKT point, recall that $v_i^+ = \frac{\partial M^+(\theta^j)}{\partial \theta_i}$ and $v_i^- = \frac{\partial M^-(\theta^j)}{\partial \theta_i}$ for $i = 1,2$ with $v_1^+ > 0$, $v_2^+ < 0$ and $v_2^- > 0$. Since $(1,0)$ lies in the convex hull of $v^+$ and $v^-$, there exists a unique $\beta$ and $v_1$ with $0 < \beta < 1$ such that $\beta \vec{v}^+ + (1-\beta)\vec{v}^- = (v_1, 0)$. Let $u_1 = 0$, $u_2 = \theta_1^j/v_1$ and $u_3 = \theta_2^j/v_1$. It is straightforward to show that these values satisfy the KKT conditions.

To complete the first case with $\beta < 1$, we need to argue that there cannot be another KKT point with $\beta = 1$. If there were another such KKT point with $\beta = 1$, then our optimization problem is a special case of the one in the previous section; hence, the solution must correspond to $\theta^b$, in which case the first component of the gradient equation would imply that $u_1 = M^-(\theta^b)/v_1$. However, this would violate the sign restriction on $u_1$ since $M^-(\theta^b)$ would be strictly positive when $\theta^j \neq \theta^b$. Consequently, when $\theta^j \neq \theta^b$, we have a unique KKT point, which must be the global minimum.

Suppose there is no solution $\theta^+$ to (4.26) on the clockwise arc following $M^+ = 0$ from $\theta^b$ to the $\theta_1$-axis. Then there cannot be a KKT point with $\beta < 1$. However, we know that there is a global minimum, which must be a KKT point. Hence, any KKT points must have $\beta = 1$. If $\beta = 1$, the optimization problem in this section becomes a special case of the optimization problem in the previous section, and it follows that there is a unique KKT point corresponding to $\theta^b$.

THEOREM 2.   *Among the jitter and bridge paths, we have*

$$(4.27) \qquad \inf_{p \in F_x}[I(p)] = \theta_1^j.$$

*Furthermore, the following conditions are equivalent:*

1. $\theta^j < \theta^b$,
2. *a jitter path is optimal in* $F_x$,
3. $\beta < 1$,
4. $M^-(\theta^b) > 0$ *and*
5. $\rho < 1$, *where*

$$\rho \equiv \frac{\bar{\lambda}_2 e^{\theta_2^j} + \mu_1 r_{1,2} e^{-\theta_1^j + \theta_2^j}}{\mu_2 r_{2,0} e^{-\theta_2^j} + \mu_2 r_{2,1} e^{\theta_1^j - \theta_2^j}}.$$

PROOF. To show (4.27), note that $\theta^j$ determines the global optimum, so we need only evaluate the objective function at the point determined by $\theta^j$. If $\theta^j = \theta^b$, the objective function becomes identical to that of the previous section, which we know evaluates to $\theta_1^j = \theta_1^b$. If $\theta^j < \theta^b$, the objective function simplifies to $\theta_1^j$.

The first three conditions are clearly equivalent from the discussion prior to the statement of the theorem. For the fourth condition, it suffices to show that the directional derivative of $M^-$ at the origin tangent to the level curve $M^+ = 0$ and going counterclockwise is strictly negative; that is, $M^- = 0$ pierces the egg $M^+ = 0$. Since $M^-(0,0) = 0$, it would follow from continuity that the level curves for $M^+ = 0$ and $M^- = 0$ would have to intersect along the segment $M^+ = 0$ somewhere between $\theta^b$ and the origin going clockwise from $\theta^b$. (From the proof of the previous theorem, we could further restrict attention to the segment lying in the first quadrant.)

If $(x^+, y^+)$ denotes the gradient of $M^+$ at the origin, then $(-y^+, x^+)$ is tangent to $M^+ = 0$ at the origin and points in a counterclockwise direction. If $(x^-, y^-)$ denotes the gradient of $M^-$ at the origin, then the directional derivative of $M^-$ at the origin tangent to the level curve $M^+ = 0$ and going counterclockwise is $-y^+ x^- + x^+ y^- = \mu_2 \bar{\lambda}_1 - \mu_1 \bar{\lambda}_1 r_{1,2} + \mu_1^* \bar{\lambda}_2 - \mu_1^* \mu_2 - \mu_1 \bar{\lambda}_2 + \mu_2 \bar{\lambda}_2 r_{2,1} + r_{1,2} \mu_1^* \bar{\lambda}_1 + r_{1,2} \mu_1^* \mu_2 r_{2,1}$. This directional derivative is a strictly decreasing function of $\mu_1^*$ since the derivative with respect to $\mu_1^*$ is $\mu_2(1 - r_{1,2} r_{2,1})(\rho_2 - 1) < 0$. If $\mu_1^* = (\lambda_1 - \rho_2 \mu_1)/(1 - \rho_2)$, then the directional derivative is zero and $M^+ = 0$ and $M^- = 0$ are tangent at the origin. However, from Proposition 1, stability requires that $\mu_1^* > (\lambda_1 - \rho_2 \mu_1)/(1 - \rho_2)$, which ensures that $M^-$ is decreasing in the direction $(-y^+, x^+)$ and pierces the egg $M^+ = 0$.

To complete the proof, using (4.6), it is straightforward to show that the fifth condition is equivalent to $v_2^+(\theta^j) < 0$; recall that $v_2^+(\theta^b) = 0$ and see Figure 2. $\qquad \square$

COROLLARY 1. *For the Jackson network with $\mu_1^* = \mu_1$, a jitter path is better than the bridge if* (2.2) *holds.*

PROOF. In the Jackson case $\exp(\theta_1^j) = \rho_1^{-1}$ and $\exp(\theta_2^j) = (r_{2,0} + r_{2,1} \rho_1^{-1})$. Using this, $\rho < 1$ by substitution and Theorem 2 gives the result. $\quad \square$

Now, we will try to derive an explicit expression for $\theta^j$ by locating the points $\theta \neq (0,0)$, where $M^+(\theta) = M^-(\theta) = 0$. Solving $M^-(\theta) = 0$ yields

$$(4.28) \qquad \exp(\theta_2) = \frac{\bar{\lambda}_1 + \bar{\lambda}_2 + \mu_1^* - \bar{\lambda}_1 e^{\theta_1} - \mu_1^* r_{1,0} e^{-\theta_1}}{\bar{\lambda}_2 + \mu_1^* r_{1,2} e^{-\theta_1}}.$$

After substituting this into $M^+(\theta) = 0$ and simplifying, we see that $x = \exp(\theta_1^j)$ must be a positive solution to the quadratic equation

$$(4.29) \qquad ax^2 + bx + c = 0,$$

where

$$
\begin{aligned}
a &= (\mu_1^* - \mu_1)(\bar{\lambda}_2 + \bar{\lambda}_1 r_{1,2})\bar{\lambda}_1 - \bar{\lambda}_2\mu_2(\bar{\lambda}_2 r_{2,1} + \bar{\lambda}_1), \\
b &= -(\mu_1^* - \mu_1)(\bar{\lambda}_2\mu_1^* + \mu_1^*\bar{\lambda}_1 r_{1,2} + \bar{\lambda}_2^2 + \bar{\lambda}_1^2 r_{1,2} + 2\bar{\lambda}_1\bar{\lambda}_2 r_{1,2} + \bar{\lambda}_1\bar{\lambda}_2 r_{1,0}) \\
&\quad + \mu_1^*\bar{\lambda}_2\mu_2(1 - 2r_{1,2}r_{2,1}) - \mu_1^*\bar{\lambda}_1\mu_2 r_{1,2}, \\
c &= (\mu_1^* - \mu_1)\mu_1^* r_{1,0}(\bar{\lambda}_2 + \bar{\lambda}_1 r_{1,2}) + (\mu_1^*)^2\mu_2 r_{1,2}(r_{1,0} + r_{1,2}r_{2,0}).
\end{aligned}
$$

Given such a solution $x$, then $\theta_1 = \log(x)$ and $\theta_2$ is determined from (4.28). Now that we have a point $\theta$ where the level curves intersect, we need to determine if it lies on the arc $M^+ = 0$ going clockwise from $\theta^b$ to the $\theta_1$-axis. Since $\theta_1 > 0$, we are only interested in $x > 1$. If, in addition, $v_2^+(\theta) < 0$, then $\theta = \theta^j$. Note that there can be at most one $x$ giving a solution that satisfies these conditions. If there is no such solution $x$, then $\theta^j = \theta^b$.

REMARK. The paper [10] provides an alternative to the Schwartz–Weiss flat boundary approach to determining the local rate function for the jitter and bridge paths. Rather than mixing the two vectors $v^+$ and $v^-$ with the weights $\beta$ and $1 - \beta$, [10] simply expresses $\Lambda^*(v, 0) = \sup_\gamma(\gamma \cdot v - \Lambda(\gamma))$, where $\Lambda(\theta) = \log(r(\hat{J}_\gamma))$, $r(\hat{J}_\gamma)$ is the spectral radius of $\hat{J}_\gamma$ and $\hat{J}_\gamma$ is the Feynmann–Kac transform of the kernel $J$ of the Markov additive process associated with the modified Jackson network as defined in [9]. If we could show that the mapping from $\gamma$ to $v$ is a smooth bijection, then since $\Lambda(\gamma)$ is convex, we could represent

$$\inf_{v \geq 0} \frac{\Lambda^*(v, 0)}{v} = \inf_\gamma\left(\gamma - \Lambda(\gamma)\Big/\frac{\partial\Lambda(\gamma)}{\partial\gamma}\right).$$

Taking derivatives, the minimum occurs when

$$\Lambda(\gamma)\frac{\partial^2\Lambda(\gamma)}{\partial\gamma^2}\Big/\left(\frac{\partial\Lambda(\gamma)}{\partial\gamma}\right)^2 = 0.$$

Since $\Lambda(\gamma)$ is convex by Lemma 2 in [10] and since $\frac{\partial\Lambda(\gamma)}{\partial\gamma} = v > 0$, it follows that the minimum occurs when $\Lambda(\gamma) = 0$ and the large deviation rate is the associated $\gamma$. This just means that the large deviation rate is the choice of $\gamma$ which sets $r(\hat{J}_\gamma) = \exp(\Lambda(\gamma)) = 1$, which agrees with the results in [9].

The problem is the assumption that the mapping from $\gamma$ to $v$ defined by $\sup_\gamma(\gamma \cdot v - \Lambda(\gamma))$ is smooth. Lemma 3.3 in [12] shows that if $\gamma \in \mathcal{U}_r$, where

$$\mathcal{U}_r = \{\gamma : \varphi(\gamma, \Lambda) = 1 \text{ for some } \Lambda = \Lambda(\gamma) < \infty\},$$

and if $(\gamma, \Lambda(\gamma)) \in \mathcal{W}$, as defined in [12], then $\Lambda(\cdot)$ is differentiable at $\gamma$. Unfortunately, we are particularly interested in cases when $\gamma$ may not belong to $\mathcal{U}_r$ because these give rise to bridges. It appears that in the two-dimensional case differentiability may follow from the explicit description of $\Lambda(\gamma)$ (private communication with Ignatiouk-Robert), but in $n$ dimensions this is far from clear. Hence, at first blush, it appears that the optimization problem would be more difficult to solve than the optimization problem arising from the Schwartz–Weiss approach where $\Lambda^*$ is represented as a convex combination of smooth pieces.

Another advantage of the Schwartz–Weiss approach is that it distinguishes between the bridge ($\beta = 1$) and the jitter ($\beta < 1$) paths. It is not clear that the approach in [10] makes this distinction. On the other hand, the results in [10] are much more general since they apply in $n$ dimensions and apply to permeable or impermeable boundaries. Note that the condition $M^+(\theta^b) > M^-(\theta^b)$ in Proposition 10 in [10] inspired condition 4 in Theorem 2.

4.3. *Cascade paths climbing the $y$-axis.* Finally, we consider fluid paths that go up the $y$-axis to a height $(0, h)$ and then go down and across $(1, 0)$. For $h > 0$, we will refer to these paths as cascade paths since the customers build up in the second queue and then "cascade" into the first queue for the large deviation. The path with $h = 0$ is none other than the bridge path. Note that a bridge path up the $y$-axis followed by a cascade is not optimal because this would give a nonlinear large deviation path in a domain with constant jump rates. Consequently, the least action path jitters up the $y$-axis and then cascades. Let $F_c \subset F$ be the set of cascade paths and the bridge path.

We wish to find conditions for a jitter path along the $y$-axis when $\rho_2 < 1$, but $\rho_1$ may be greater than one. The investigation of jitter paths in Section 4.2 can be used provided we interchange the $x$- and $y$-axes. Let $\tilde{\theta} = (\tilde{\theta}_1, \tilde{\theta}_2)$ be the solution analogous to $\theta^j$ in Theorem 2. Thus, equivalent to (4.26), we have

$$M^+(\tilde{\theta}_1, \tilde{\theta}_2) = \tilde{M}(\tilde{\theta}_1, \tilde{\theta}_2) = 0,$$

where

$$\tilde{M}(\tilde{\theta}_1, \tilde{\theta}_2) := \bar{\lambda}_1(e^{\tilde{\theta}_1} - 1) + \bar{\lambda}_2(e^{\tilde{\theta}_2} - 1) + \mu_2 r_{2,0}(e^{-\tilde{\theta}_2} - 1) + \mu_2 r_{2,1}(e^{\tilde{\theta}_1 - \tilde{\theta}_2} - 1).$$

Subtracting $\tilde{M}$ from $M^+$ gives $\exp(\tilde{\theta}_1) = r_{1,0} + r_{1,2} \exp(\tilde{\theta}_2)$. Substituting this into $\tilde{M}(\tilde{\theta}) = 0$ gives $\exp(\tilde{\theta}_2)$ equal to 1 or $\mu_2/\lambda_2$, and the former solution can be eliminated since it corresponds to $\tilde{\theta}_1 = \tilde{\theta}_2 = 0$. In order to have a jitter path, condition 4 in Theorem 2 must hold; that is,

$$\tilde{\rho} := \frac{\bar{\lambda}_1 e^{\tilde{\theta}_1} + \mu_2 r_{2,1} e^{-\tilde{\theta}_2 + \tilde{\theta}_1}}{\mu_1 r_{1,0} e^{-\tilde{\theta}_1} + \mu_1 r_{1,2} e^{\tilde{\theta}_2 - \tilde{\theta}_1}} < 1.$$

Substituting $\exp(\tilde{\theta}_2) = \mu_2/\lambda_2$ gives $\tilde{\rho} = \lambda_1 \exp(\tilde{\theta}_1)/\mu_1 < 1$; that is, (2.3) holds.

We conclude there is no cascade path unless (2.3) holds, but if it does, then the large deviations rate of paths from $(0,0)$ to $(0,h)$ is $h \log(\rho_2^{-1})$, which can be seen by simplifying (4.29) when $\mu_1 = \mu_1^*$ to obtain $\theta_1^j = \rho_1^{-1}$. Moreover, the large deviation rate of paths from $(0,h)$ to $(1,0)$ is

$$\inf_{v>0} \frac{\Lambda^+(v,-vh)}{v}.$$

After $(1,0)$, the process follows the natural drift path back to zero. Thus, if path $p \in F_c$ reaches height $h$, we have $I(p) = h \log(\rho_2^{-1}) + \inf_{v>0} \frac{\Lambda^+(v,-vh)}{v}$. We wish to find $\inf_{p \in F_c} I(p)$ or, equivalently, $\inf_{h \geq 0, v > 0} f(h,v)$, where $f(h,v) = h \log(\rho_2^{-1}) + \frac{\Lambda^+(v,-vh)}{v}$. Since $f$ is continuous, positive and diverges as $v \downarrow 0$ or as $v \to \infty$ or as $h \downarrow 0$, we know $\inf_{h \geq 0, v > 0} f(h,v)$ must be a local minimum of

$$(4.30) \qquad\qquad \min f(h,v)$$

$$(4.31) \qquad\qquad \text{s.t.} \quad h \geq 0.$$

This constrained nonlinear optimization problem is differentiable with a linear constraint so every local minimum must be a KKT point. In order for $(h,v)$ to be a KKT point for this problem, there needs to be a corresponding Lagrange multiplier $u$ so that $h, v$ and $u$ satisfy the following:

$(4.32) \quad uh = 0 \qquad$ (complementary slackness),

$(4.33) \quad u \geq 0 \qquad$ (sign restriction),

$(4.34) \quad \begin{pmatrix} u \\ 0 \end{pmatrix} = \nabla f(h,v) = \begin{pmatrix} \log(\rho_2^{-1}) - \theta_2^+(v,-vh) \\ M^+(\theta^+(v,-vh)) \end{pmatrix} \qquad$ (gradient equation),

$(4.35) \quad h \geq 0 \qquad$ (constraint).

As in the previous sections, there will be a unique KKT point $(h,v) \in [0,\infty) \times (0,\infty)$.

THEOREM 3. *Define $\theta^c = \theta^+(v,-vh)$, where $(v,h)$ minimizes* (4.30) *subject to* (4.31). *It follows that*

$$(4.36) \qquad\qquad \inf_{p \in F_c}[I(p)] = \theta_1^c.$$

*If $\rho_1^{-1} > r_{1,0} + r_{1,2}\rho_2^{-1}$ and $\log(\rho_2^{-1}) < \theta_2^b$, then the minimum action path in $F_c$ is a cascade of height $h > 0$ with $\theta^c = (\log(\rho_1^{-1}), \log(\rho_2^{-1}))$; otherwise, the minimum action path in $F_c$ is a bridge with $\theta^c = \theta^b$.*

PROOF. We will argue that $\theta^c$ corresponds to the unique KKT point in $(h,v) \in [0,\infty) \times (0,\infty)$; hence, this point must correspond to the global minimum. Note that the second component of the gradient equation implies that we must be looking for a solution corresponding to a point $\theta^+$ with $M^+(\theta^+) = 0$. This point $\theta^+$ determines $v$, $h$ and $u$ since $v = v_1^+(\theta^+)$, $h = -v_2^+(\theta^+)/v$, and $u = \log(\rho_2^{-1}) - \theta_2^+$.

In order to have a KKT point $\theta^+$ with $u > 0$ we would have to have $h = 0$, which corresponds to a bridge path. Hence, $\theta^+ = \theta^b$ and by (4.34), this occurs if and only if (2.3) is false or $\log(\rho_2^{-1}) > \theta_2^b$. (4.36) follows by substitution.

Alternatively, in order to have a solution with $u = 0$, we would have to have $h > 0$, which corresponds to a cascade path. Hence, $\theta^+ = \theta^c$ and this occurs if and only if (2.3) is true and $\log(\rho_2^{-1}) \le \theta_2^b$. Next, (4.34) requires $\theta_2^+ = \log(\rho_2^{-1}) \le \theta_2^b$ and $M^+(\theta^+(v, -vh)) = 0$. Moreover, $\nabla M^+(\theta^+) = (v, -vh)$. The two solutions to $M^+(\theta^+(v, -vh)) = 0$ and $\theta_2^+ = \log(\rho_2^{-1})$ are $(\log(\rho_1^{-1}), \log(\rho_2^{-1}))$ and $(\log(r_{1,0} + r_{1,2}\rho_2^{-1}), \log(\rho_2^{-1}))$.

Substituting the second solution into (4.5) and (4.6) gives

$$\nabla M^+(\theta^+) = ((\bar{\lambda}_1 + \lambda_2 r_{2,1})e^{\theta_1^+} - \mu_1 e^{-\theta_1^+}(r_{1,0} + \rho_2^{-1}r_{1,2}),$$

$$\bar{\lambda}_2\rho_2^{-1} - \lambda_2 r_{2,0} + \mu_1 r_{1,2}\rho_2^{-1}e^{-\theta_1^+} - \lambda_2 r_{2,1}e^{\theta_1^+}).$$

Notice that the first coordinate equals $\lambda_1 e^{\theta_1^+} - \mu_1 = \lambda_1((r_{1,0} + r_{1,2}\rho_2^{-1}) - \rho_1^{-1})$. This is negative if (2.3) is true and this is impossible since $\nabla M^+(\theta^+) = (v, -vh)$ with $v > 0$.

Substituting the first solution into (4.5) and (4.6) gives

$$\nabla M^+(\theta^+) = (\mu_1 - \lambda_1 r_{1,0} - \lambda_1 r_{1,2}\rho_2^{-1}, \mu_2 - \lambda_2 r_{2,0} - \lambda_2 r_{2,1}\rho_1^{-1}) = (v, -vh).$$

The first component is $\lambda_1(\rho_1^{-1} - (r_{1,0} + r_{1,2}\rho_2^{-1}))$ and this is positive since (2.3) is true. The second component of $\nabla M^+(\theta^+)$ is negative since $\theta^+$ is on the level curve $M^+(\theta) = 0$ and $\theta_2^+ < \theta_2^b$. (2.2) fails. Equation (4.36) follows by substitution.

<div align="right">□</div>

4.4. *Summary of flat boundary approach.* Combining yields the following result.

THEOREM 4. *For the modified Jackson network,*

$$\lim_{\ell \to \infty} \frac{1}{\ell} \log P(W \in F(\ell)|W(0) = (0,0))$$

$$= \lim_{\ell \to \infty} \frac{1}{\ell} \log P(W_\ell \in F(1)|W_\ell(0) = (0,0))$$

$$= -\inf_{p \in F}[I(p)] = -\inf_{p \in F_i \cup F_x \cup F_c}[I(p)].$$

*If $\theta_2^j < \min\{\log(\rho_2^{-1}), \theta_2^b\}$, then minimum action is $\theta_1^j$ and the minimal action path is a jitter path along the x-axis. If $\log(\rho_2^{-1}) < \min\{\theta_2^j, \theta_2^b\}$, then the minimal action is $\theta_1^c$ and the minimal action path is a cascade path that initially climbs the y-axis. Otherwise, $\theta_1^b = \theta_1^j = \theta_1^c$, the minimal action is $\theta_1^b$, and the minimal action path is a bridge.*

PROOF.  The only part of the proof that is not straightforward is showing that if $\log(\rho_2^{-1}) < \min\{\theta_2^j, \theta_2^b\}$, then the condition $\rho_1^{-1} > r_{1,0} + r_{1,2}\rho_2^{-1}$ of Theorem 3 holds automatically. To prove this, we will show that if $\rho_1^{-1} \leq r_{1,0} + r_{1,2}\rho_2^{-1}$ and $\log(\rho_2^{-1}) < \theta_2^b$, then $\theta_2^j < \log(\rho_2^{-1})$. In other words, a jitter path will be the minimal action path—not a cascade.

The (convex) function $g^-$ defined by $M^-(\theta_1, g^-(\theta_1)) = 0$ is given at (4.28). Simplifying,

$$g^-(\theta_1) = 1 + \frac{\mu_1^*(1 - e^{-\theta_1}) + \bar\lambda_1(1 - e^{\theta_1})}{\bar\lambda_2 + \mu_1^* r_{1,2} e^{-\theta_1}}.$$

Differentiating $g^-$ with respect to $\mu_1^*$ gives $(1 - \exp(-\theta_1))/(\bar\lambda_2 + \mu_1^* r_{1,2} \times \exp(-\theta_1))^2$ and this is strictly positive for $\theta_1 \in [0, \log(\mu_1^*/\bar\lambda_1)]$, where the endpoints are the zeros of $g^-$. Thus, $g^-$ is strictly increasing in $\mu_1^*$ on $[0, \log(\mu_1^*/\bar\lambda_1)]$. Furthermore, $\lim_{\mu_1^* \to \infty} g^-(\log(r_{1,0} + r_{1,2}\rho_2^{-1})) = \log(\rho_2^{-1})$, which is a point on $M^+ = 0$. If $\rho_1^{-1} \leq r_{1,0} + r_{1,2}\rho_2^{-1}$, then the easternmost solution $\theta$ to $M^+(\theta) = 0$ with $\theta_2 = \log(\rho_2^{-1})$ is $(\log(r_{1,0} + r_{1,2}\rho_2^{-1}), \log(\rho_2^{-1}))$. It follows that $M^+ = 0$ and $M^- = 0$ intersect at a point in the positive quadrant going clockwise from $(\log(r_{1,0} + r_{1,2}\rho_2^{-1}), \log(\rho_2^{-1}))$ along $M^+ = 0$ before hitting the axis. If, in addition, $\log(\rho_2^{-1}) < \theta_2^b$, then the point of intersection defines $\theta^j$, which must be the minimal action path.

There is one boundary case which also must be eliminated. We must show it is impossible that $\theta^j = \theta^c < \theta^b$. If this did happen, then $\theta^j = \theta^c = (\log(\rho_1^{-1}), \log(\rho_2^{-1}))$. Since $\theta^j$ lies on $M^+$ and $M^-$, it follows that $g^-(\log(\rho_1^{-1})) = \rho_2^{-1}$. Since $\theta^j < \theta^b$, it follows that $\rho < 1$ and this reduces to $\rho_2^{-1} < r_{2,0} + r_{2,1}\rho_1^{-1}$. Now notice that when $\mu_1^* = \mu_1$, that is, in the Jackson case, we have $g^-(\log(\rho_1^{-1})) = r_{2,0} + r_{2,1}\rho_1^{-1}$ and this is strictly greater than $\rho_2^{-1}$ as we just showed. This eliminates this case if $\mu_1^* = \mu_1$. Moreover, we showed above that $g^-(\theta_1)$ is increasing in $\mu_1^*$, so it will never be possible to solve $g^-(\log(\rho_1^{-1})) = \rho_2^{-1}$.  □

COROLLARY 2.  *For the Jackson network with $\mu_1^* = \mu_1$, if (2.2) fails and $(r_{2,0} + r_{2,1}\rho_1^{-1}) > \rho_2^{-1}$, then a cascade path is optimal, whereas if (2.2) fails and $(r_{2,0} + r_{2,1}\rho_1^{-1}) = \rho_2^{-1}$, then a bridge path is optimal. In both cases the minimal action is $\log(\rho_1^{-1})$. (This is no surprise since the steady state $\pi$ is of product form even in the cascade or bridge cases.)*

Proof. In the Jackson case $\exp(\theta_1^j) = \rho_1^{-1}$ and $\exp(\theta_2^j) = (r_{2,0} + r_{2,1}\rho_1^{-1})$. Hence, if $(r_{2,0} + r_{2,1}\rho_1^{-1}) > \rho_2^{-1}$, then $\exp(\theta_2^j) > \rho_2^{-1}$. Using this, $\rho \geq 1$ by substitution and Theorem 2 shows $\theta^j = \theta^b$. This means $\exp(\theta_2^b) > \rho_2^{-1}$ and, moreover, if (2.2) fails, then (2.3) must hold. By Theorem 3 it follows the minimum action path is a cascade and $\theta_1^c = \log(\rho_1^{-1})$.

If $(r_{2,0} + r_{2,1}\rho_1^{-1}) = \rho_2^{-1}$, then by the above argument $\exp(\theta_2^j) = \rho_2^{-1}$. Hence, by Theorem 3, the bridge path is optimal and $\theta_1^b = \theta_1^j = \log(\rho_1^{-1})$. □

**5. Another model exhibiting the bridge phenomenon.** We can use the above methods to obtain the rough asymptotics of $\pi(\ell, y)$ for a fork network as introduced by Flatto and Hahn [7]. This model was later extended and described as the bathroom problem in Chapter 16 of [15]. Couples arrive at a cinema according to a Poisson process with rate $\nu$ and immediately visit the men's and ladies' room. The service rate the men's queue $\alpha$, while the rate at at the lady's queue is $\beta$. There are also separate arrival streams, with rate $\eta$ for single women and rate $\lambda$ for single men.

We are interested in a large deviation of the men's queue, so let this be the first queue and the ladies' queue the second. We could use the above methods to obtain the rough asymptotics of $\pi(\ell, y)$, but we don't have to because the work is already done in Section 16.2 in [15]. In the (unrealistic) case $\eta + \alpha\nu/(\lambda + \nu) < \beta$, they show the most likely path for the men's queue to reach a high level $\ell$ before returning to zero is a jitter path along the $x$-axis and the large deviation rate is $\alpha/(\lambda + \nu)$ in agreement with [11]. In the (more realistic) case when $\eta + \alpha\nu/(\lambda + \nu) > \beta$, the results in [15] show the most likely path for the men's queue to reach a high level $\ell$ before returning to zero is a path through the interior (again in agreement with [11]). Moreover, among jitter paths, the large deviation rate (for the more realistic case) is minimized by paths spending zero time on the boundary.

This just means that if, in addition to the men's queue getting big, we require that the ladies' queue remain small, then the large deviation path is a bridge path. The exact asymptotics of $\pi(\ell, y)$ are calculated in [9] and are found to agree with those in [7]. The derivation of these sharp asymptotics confirms our intuition about the bridge behavior.

**Acknowledgments.** David R. McDonald wishes to thank the members of the Department of Mathematics of the Indian Institute of Science in Bangalore for their hospitality while working on this paper.

SCHOOL OF INDUSTRIAL
  AND SYSTEMS ENGINEERING
GEORGIA INSTITUTE OF TECHNOLOGY
765 FERST DRIVE
ATLANTA, GEORGIA 30332-0205
USA

DEPARTMENT OF MATHEMATICS
  AND STATISTICS
UNIVERSITY OF OTTAWA
P.O. BOX 450
OTTAWA, ONTARIO K1N6N5
CANADA
E-MAIL: dmdsg@uottawa.ca